\documentclass[amsfonts,amssymb]{article}

\usepackage{amstext,amsfonts,amsmath,amssymb}

\begin{document}
\title{\bf Probability and strategy in a variation of a blackjack game}
\date{}

\author{
Zoran Stani\' c
\\\footnotesize{Faculty of Mathematics, University of Belgrade}\\
\footnotesize{Studentski trg 16, 11000 Belgrade,
Serbia}\\\footnotesize{ {\it Emails:} zstanic@matf.bg.ac.rs;
zstanic@math.rs}}

\date{}

\maketitle

\begin{abstract}\noindent
The subject of this paper is a variation of a blackjack game,
mainly popular in some parts of Europe where it is known as einz
(in German slang: one). We describe the rules of this game,
indicate its main characteristics, give some probabilities,
suggest the strategy, and consider some typical situations.

\end{abstract}

\section{Introduction}

Einz is a game for two or more players played with one or more
decks of 52 cards. Face cards king, queen, and jack  are counted
as 4, 3, and 2 points, respectively. An ace is counted as 11
points, and the remaining cards are counted as the numeric value
shown on the card. The value of a hand is simply the sum of the
point counts of each card in the hand. Any hand with value of 21
is called an {\it einz}, as well as the hand with exactly two
cards and value of 22 (these cards must be two aces). The object
of the game is to reach an einz or to reach a score higher than
the opponents without exceeding 21. Scoring higher than einz (a
{\it busting}) results in a loss.

There are two basic types of this game:

\begin{description}
    \item[{\rm (1)}] game with no dealer (or {\it open game}) in which
    a player is playing against all other players, and
    \item[{\rm (2)}] game with a dealer
(or {\it dealer game}) in which a player is playing against the
dealer, but not against any other players.
\end{description}

In an open game, the cards are dealt by a different player in each
round. The player who is dealing the cards does not have any
special status in relation to the others. At begin, every player
is dealt an initial hand of two cards. Then, each player, one by
one, is getting the additional cards in order to achieve his total
score. All cards in the hand are visible to himself only, and a
player can stand at any time. A player must announce a bust
immediately after it happens, and in this case he lose
automatically. If a player gets an einz, he is showing his hand,
and he wins the round (no matter if any other players did not got
a chance to make their total scores). If none player got an einz,
the winner is whoever has closest to a total of 21.

There is an additional rule usually known as {\it changing on 14}:
each player can start over with getting the new cards whenever his
total score equals 14. If he decide to get the new cards, he is
showing the hand, and continue with the new one.

There are few variations of a dealer game, and they will be
explained and discussed in Section 3.

\medskip

There are many, more or less popular, variations of the blackjack
game. The main characteristics of the einz game are the distinct
values of face cards, immediate announcing of a bust or an einz,
the einz with 2 aces, and changing on 14. Note that some of these
rules appear in some other variations of the game. On the other
hand, a very common rules of many popular blackjack games
(splitting the hand, soft and hard ace, etc -- see \cite{Jam,
Wat}) are not present here.

\medskip

The strategy in usual blackjack and its simple variations (played
in worldwide casinos) is well developed. It includes basic
probabilistic methods, as well as so called card counting
\cite{Jam, Wat, Ols} - a card game strategy used to determine
whether the next hand is likely to give a probable advantage to
the player or to the dealer. The most common variations of card
counting are based on statistical evidence that high cards,
especially aces and all cards counted as 10 points (there, face
cards are all counted as 10 points, as well as 10s), benefit the
player more than the dealer, while the low cards, especially 4s,
5s, and 6s, help the dealer while hurting the player. Thus, a high
concentration of high cards in the deck increases the player's
chance, while low cards benefit the dealer. Contrary to the usual
blackjack, in einz we have more low than high cards, and the usual
counting principles cannot be transferred. But in any case, here
we consider the einz as an easy relaxing game, and therefore
through the whole paper we will not include the counting strategy
nor playing for the money. In order to avoid the possibility of
counting, we will assume that each round starts with a full set of
cards.

\medskip

In the following two sections we consider the open and the dealer
game, respectively. Some typical situations are discussed in the
last section. We use the usual terminology, and techniques in
probability that can be found in any complete book concerning this
discipline, say \cite{prob} or \cite{prob2}.

\section{Open game}

\subsection{When to stand?}

Here we consider the question: when the player should stand in the
open einz game? According to the rules, two obvious candidates for
the standing strategy arise

\begin{description}
    \item[{\rm (1)}] the player should stand when its score is
    at least 17, or
    \item[{\rm (2)}] the player should stand when its score is
    at least 18.
\end{description}

If the player uses the standing strategy (1) in every round, i.e.
if he is taking the cards until his hand busts or achieves a score
of 17 or higher in every round, we will say that "he always stands
on 17". Similarly, if he uses the strategy (2), we will say that
"he always stands on 18".

We compute the probabilities of all possible outcomes in both
suggested strategies. To simplify the computation, we assume that
the player never changes on 14. An einz in two cards can be
achieved if these two cards are either both aces or an ace and a
ten. Therefore, the probability that the player gets an einz in
two cards is easily computed and it is equal to $0.016$. The
remaining situations are resolved in the similar way (we use the
computer to complete the whole computation procedure). The results
are given in Table 1 and Table 2.

\begin{center}

\begin{tabular}{|c|c|c|c|c|c|}
  \hline
  $\sharp$ cards / score & 17 & 18 & 19 & 20 & einz \\
  \hline
  2 & 0.036 & 0.027 & 0.024 & 0.016 & 0.016 \\
  3 & 0.074 & 0.067 & 0.051 & 0.035 & 0.038 \\
  4 & 0.040 & 0.042 & 0.056 & 0.036 & 0.027 \\
  5 & 0.010 & 0.009 & 0.013 & 0.012 & 0.008 \\
  $>5$ & $<0.001$ & $<0.001$ & $<0.002$ & $<0.002$ & $<0.002$ \\
  any & 0.161 & 0.147 & 0.145 & 0.100 & 0.092 \\
  \hline
\end{tabular}

\bigskip

{\footnotesize Table 1: Probabilities of scores between 17 and
einz for the player who always stands on 17.}

\bigskip

\begin{tabular}{|c|c|c|c|c|c|}
  \hline
  $\sharp$ cards / score  & 18 & 19 & 20 & einz \\
  \hline
  2  & 0.027 & 0.024 & 0.016 & 0.016 \\
  3  & 0.067 & 0.056 & 0.040 & 0.044 \\
  4  & 0.042 & 0.066 & 0.046 & 0.038 \\
  5  & 0.009 & 0.018 & 0.017 & 0.014 \\
  $>5$ & $<0.001$ & $<0.003$ & $<0.003$ & $<0.003$ \\
  any &  0.147 & 0.167 & 0.122 & 0.114 \\
  \hline
\end{tabular}

\bigskip

{\footnotesize Table 2: Probabilities of scores between 18 and
einz for the player who always stands on 18.}

\end{center}

By summation of the values in the last row of Table 1, we get that
the probability of score 17 or better in a single hand is equal to
0.645, and thus the probability of busting is 0.355. Considering
Table 2, we get that these two probabilities are in this case
equal to $0.550$ and $0.450$, respectively. In the first look the
first strategy seems to be better, but let us consider both of
them in a real game with two players.

Assume first that both players always stand on 17. Then, according
to the rules, the first player wins in the following situations:
(a) if he gets an einz, (b) if his score is between 17 and 20, and
the second player busts, (c) if his score is between 17 and 20 and
it is higher than score of the second player. Using Table 1, we
get the probabilities of all these outcomes: (a) 0.092, (b) 0.196,
and (c) 0.114, which together give the first player's winning
probability equal to 0.402. The round is tied if both players have
equal score between 17 and 20, the probability of that is 0.079.
Finally the second player's winning probability is equal to
1-0.402-0.079=0.520. These results are given in the first row of
Table 3.  In the same way we get the remaining data of the same
table.

\begin{center}

\begin{tabular}{|c|c|c|c|c|c|}
  \hline
  result /  stand on  & 17 vs 17 & 17 vs 18 & 18 vs 17 & 18 vs 18\\
  \hline

  player 1 wins & 0.402 & 0.379 &  0.399 & 0.373\\
  tied & 0.079  & 0.058 & 0.058&0.064  \\
  player 2 wins & 0.520  & 0.563 & 0.543& 0.562  \\

  \hline
\end{tabular}

\bigskip

{\footnotesize Table 3: Probabilities in the game with two players
with given strategies.}

\end{center}

Regarding Table 3, the strategy (2), i.e. "always stand on 18", is
clearly better choice for the second player, while for the first
player the strategy (1) has a tiny advantage.

The game with more than 2 players can be considered in the similar
way. Here we give just one result: in the game with 3 players such
that the first two players always stand on 17, and the third
player always stand on 18, the first player's winning probability
is 0.1966, the second's is 0.1881, and the third's is 0.3494,
while at least two players have mutually  equal winning score with
probability of 0.2659. In this game the first player has a bit
better chances than the second one (that is because the round is
over if the first player gets an einz, while if he busts, the
second player can bust as well with equal probability). Note also
that the third player has a clear advantage, and the same holds in
any game with $k~(k>3)$ players: the $k$-th player has advantage
against any other. Observe finally that the game with 3 players
reduces to a game with 2 players (along with the probabilities
presented above) whenever the first player busts.

\bigskip

\noindent {\it Conclusion: The first player should stand on  17,
the second player should stand on  18. In the game with more than
two players, the last player has a clear advantage.}

\subsection{Change on 14 or not?}

We consider the question in the title of this part. In order to
get the exact probabilities we consider an einz game with 52
cards, and we also assume that the player gets 10 and 4 in his
first two cards (the conclusion given in the end is the same in
any other case, i.e. it does not depend on these assumptions).

Assume first that the player always stands on 17. Then he does not
bust if his next card is counted as 3, 4, 5, 6, or 7 (recall that
queens are counted as 3, as well as 3s), or if his next two cards
are counted as 2 and $k~(k=2,3,4,5)$ and if they are taken in the
given order.

We compute the probabilities: $P(3)=\frac{8}{50},
P(4)=\frac{7}{50}, P(5)=P(6)=P(7)=\frac{4}{50}$, and
$P(2,2)=\frac{56}{2500}, P(2,3)=\frac{64}{2500},
P(2,4)=\frac{56}{2500}, P(2,5)=\frac{32}{2500}$. Altogether, the
probability of 17 or better score is equal to the sum of these
 probabilities and it is equal to 0.624, i.e. it is less than the
 probability of the same result if he changes on 14 (this probability is slightly different from 0.645 (this results is obtained in the previous subsection) since two cards, 10 and 4,
 are removed from the deck, but it is very close to this number).
 Therefore, change on 14 is a bit better choice for this player.

\medskip

Assume now that the player always stands on 18. If he decides to
get another one or two cards, we get that the probability of 18 or
better score is equal to 0.514, i.e. it is again less than the
probability of the same score if he changes on 14.

\bigskip

\noindent {\it Conclusion: Regardless of his standing strategy, a
player will tiny increase his chances if he changes on 14.}

\subsection{Probability and mathematical expectation of possible scores}

What we can say about his hand if  a player stands after he took a
few cards? Clearly, he did not busted and he did not got an einz
(since both results must be announced), and so he has some high
score. It is an intriguing question: which score he probably has?
In other words, it would be interesting to compute the
probabilities of all possible scores. This computation depends on
(earlier discussed) playes's  standing strategy, and the number of
cards that he took.

Assume that the player always stands on 17, and that he stands
after he took 2 cards. Then, according to the first row of Table
1, the probability that he has 17 is equal to
$0.036/(0.036+0.027+0.024+0.016)=0.350$. In the similar way, we
compute the remaining probabilities (only if he took between 2 and
5 cards, the remaining situations are very rear), and the results
are given in the first part of Table 4. In the second part of this
table we give the results for the player who always stands on  18.

\begin{center}

\begin{tabular}{|c|c|c|c|c|c|c|}
  \hline
  $\sharp$ cards / score & 17 & 18 & 19 & 20 \\
  \hline
  2 & 0.350 & 0.262 & 0.233 & 0.156  \\
  3 & 0.326 & 0.295 & 0.225 & 0.154  \\
  4 & 0.230 & 0.241 & 0.322 & 0.207  \\
  5 & 0.227 & 0.205 & 0.295 & 0.273  \\
  \hline\hline
2 &  & 0.403 & 0.358 & 0.239  \\
  3 &  & 0.411 & 0.344 & 0.245  \\
  4 &  & 0.273 & 0.429 & 0.299  \\
  5 & & 0.205 & 0.409 & 0.386  \\

  \hline
\end{tabular}

\bigskip

{\footnotesize Table 4: Probability of score if the player stands
with $k~(k=2,\ldots,5)$ cards in the hand.}

\end{center}

As we can see, if the player took only two cards, there are small
chances that he has the highest possible score (i.e. 20), but
these chances increase with the number of cards. The remaining
possible scores can be analyzed in the similar way. Another
application of this computation is given below.

Assume, for example, that we have a game with two players, and let
say that we know that both of them always stand on 17 (different
possibilities can be considered in the similar way). If both of
them stand after they took equal number of cards then their
chances to win are also equal (nobody got an einz and nobody
busted). What if they took a different number of cards? The
probabilities of the possible results are given in Table 5. There,
regardless of the number of cards, the probability of tied round
is about 0.25 (i.e. 1 of 4 rounds should be tied). On the other
hand there is a large difference in wining probabilities.

\begin{center}

\begin{tabular}{|c|c|c|c|c|c|c|}
  \hline
  result / $\sharp$ cards & 2 vs 3 & 2 vs 4 & 2 vs 5 & 3 vs 4 & 3 vs 5 & 4 vs 5 \\
  \hline
  player 1 wins & 0.361 & 0.293 & 0.273 & 0.296 & 0.276&0.344 \\
  tied & 0.268 & 0.251 & 0.244 & 0.250 &0.243& 0.253 \\
  player 2 wins & 0.371 & 0.456 & 0.483 & 0.454 &0.481& 0.403\\

  \hline
\end{tabular}

\bigskip

{\footnotesize Table 5: Probabilities in the game with two players
(both always stand on 17)  if they achieved their scores between
17 and 20, the first with  $k$, and the second with $l$ cards
($2\leq k<l\leq 5$).}

\end{center}

The mathematical expectation (if necessary, see definition in
\cite[Chapter 6]{prob}), or in the case of this game, the expected
average score of the player who stands after he took 2 cards (if
he always stands on 17) is (compare the first row of Table 4)
$E(2;17-20)=0.350\cdot 17+0.262\cdot 18+0.230\cdot 19+0.227\cdot
20=18.156$. The remaining mathematical expectations (computed
including or not including an einz into the high score) in the
case of two standing strategies are given in Table 6.

\begin{center}

\begin{tabular}{|c|c|c|c|c|c|c|}
  \hline
  $E$ / $\sharp$ cards & 2  & 3 & 4 &  5 & any\\
  \hline
  17-20 & 18.156 & 18.207 & 18.506 & 18.614& 18.332 \\
  17-einz & 18.571 & 18.608 & 18.841 & 18.981& 18.707 \\
\hline\hline
18-20 & 18.836 & 18.834 & 19.026 & 19.182& 18.943 \\
  18-einz & 19.256 & 19.295 & 19.417 & 19.621& 19.369 \\

  \hline
\end{tabular}

\bigskip

{\footnotesize Table 6: Average value of the high score achieved
with $k~(2\leq k\leq 5)$ cards.}

\end{center}

\bigskip

\noindent {\it Conclusion: When a player stands the following rule
can be applied: more cards in his hand - better result can be
expected. The average score also increases in the number of
cards.}

\section{Dealer game}

In this type of the game the dealer is dealing the cards in every
round, and therefore his hand is resolved after all players have
finished playing. In many variations of the blackjack game the
dealer's playing is strictly determined  by certain rules
directing his every move in the game. Commonly, the dealer game in
einz is less popular than the open game, but here we present a
several possibilities for the dealer's playing describing at the
same time the whole game.

First, it is usual that the additional possibility, changing on
14, is not allowed to the dealer (the same holds for the similar
rules in better known variations of this game (like "splitting the
hand" \cite{Jam, Wat})). It is also usual that all cards dealt are
visible to everyone.

The possibilities of the remaining rules follow:

\begin{description}
    \item[{\rm (i)}]~ Here, contrary to the open game the busting and getting an einz
    do not mean automatic lose or win. In this case, if the dealer
    equals the player's result (both busted, or got the equal score including the einz) the game is tied. The dealer always stands on, say, 17.
 This variation is very similar to
many of popular types of the blackjack game. But, without
possibility that anyone changes on 14, these rules give the equal
chances to a player and dealer whenever the player applies the
dealer's standing strategy. If the player is allowed to use the
opportunity to change on 14, then his winning probability is a bit
higher, which is not usual, and would be corrected by some
additional rule concerning the paying the wins (this part of the
game is beyond the subject of this paper).

    \item[{ \rm(ii)}] ~Regarding a player, all rules of the open game remain unchanged. The dealer always stand on 17.
    By the rules of this variation, the dealer wins if (a) the player busts, (b) the
dealer's score is strictly higher than the player's score.

According to these rules, none round is tied. If the player always
stands on 17, the probability that he busts is 0.355, and if he
does not bust the probability that the dealer gets a higher score
is 0.165. Thus, the player's  winning probability is 0.480
(against the dealer's 0.52) -- see Table 3. Using the opportunity
to change on 14, and by changing his standing strategy, the player
can increase his chances (but not over 0.5). Since, it is usual
that the dealer has a better chances, but close to 0.5, this
variation of dealer game is the author's recommendation to be
used.

If the player always stands on at least 18, his winning
probability is 0.458 (less than the above).

    \item[{ \rm (iii)}] The same as (ii), with exception that the dealer
    always stands on 17.

    In this case the dealer's chances are better than the above. Regardless of player's standing strategy, the  dealer's
    winning probability is about 0.563.
\end{description}

\section{Typical situations}

There are many typical situations in both open and dealer game
that can be resolved using results of Section 2. We consider just
few of them.

\begin{enumerate}
    \item {\it Open game with 8 card decks. Three players are playing, the first two both stand with 2 cards in the hand, while the third scored 18 in his 2 cards. What should he do?}

Every player took 2 cards, and the value of each of these 6 cards
is higher than 5 (if we naturally assume that none player stands
before 17). Assume first that the first two players both always
stand on 17. If the third player stands then, according to the
simple computation based on Table 1, his chances are about these
numbers:

\begin{itemize}
    \item win: 0.123
    \item win together with one player: 0.104
    \item all 3 players tied: 0.069
    \item lose: 0.705
\end{itemize}

If he take another card, he busts (and lose) if the value of this
card is higher than 3. He wins if the value is 3, or if the value
is 2 and the score of both previous players is less than 20. In
the remaining situations he is tied with at least one player. His
chances are:

\begin{itemize}
    \item win: 0.202
    \item win together with one player: 0.586
    \item all 3 players tied: 0.515
    \item lose: 0.688
\end{itemize}

So, in the last case the chances are better. If any of the
previous players does not stand on 17, the difference between
these chances is even bigger.

Conclusion: Although he scored 18, he should take another card!

\item {\it Open game with 8 card decks. Two players are playing,
and the first of them has 10 and 6. What should he do?}

Assume that the second player never stands before 18. If the first
player stands his chances are:

\begin{itemize}
    \item win: 0.45
    \item lose: 0.55
\end{itemize}

If he take another card we get the following probabilities:

\begin{itemize}
    \item win: 0.357
    \item tied: 0.067
    \item lose: 0.576
\end{itemize}

So, the first player has better chances if he stands on 16! But,
assume now that the second player always stands on  17. In the
similar way we get the following results. If the first player
stands, his chances are:

\begin{itemize}
    \item win: 0.355
    \item lose: 0.645
\end{itemize}

If the first player take another card his chances are:

\begin{itemize}
    \item win: 0.385
    \item tied: 0.061
    \item lose: 0.554
\end{itemize}

So in this case he has better chances if he take another card.

Conclusion: The answer depends on the standing strategy of the
second player. We conclude earlier that "stand on 18" is a better
strategy for him. But this situation shows that this strategy is
not always the best choice. We remind that the second player is
choosing his strategy when he does not know whether the first
player has the same strategy in every round or not. An additional
conclusion would be that the standing strategy should be changed
through the game.

    \item {\it Dealer game with 8 card decks. The player has 9 and 8. Discuss what should he do if the game is playing under the rules described in (ii) or (iii) (see the previous section).}

In (ii), if he stands his winning probability is 0.516 (better
chances than the dealer), and if he take another card his winning
probability is 0.437.

In (iii), if he stands his winning probability is 0.45, and if he
take another card his winning probability is 0.427.

Conclusion: In both situations he should stand.

   \item {\it Dealer game described in (ii). The player and  dealer stand, but the player has less cards in the hand. Compute the probability of the player's win.}

If the player always stands on 17 then, according to Table 5, he
has over than 50 percent chances to win. For example, if he took
2, and the dealer took 3 cards, then the player's winning
probability is 0.629, etc. The other player's standing strategy is
left to the reader.

\end{enumerate}

\end{document}